\documentclass[twoside,a4paper,12pt,centertags,reqno]{amsart} 
\usepackage{amsmath,amssymb,verbatim,vmargin}
\usepackage{color}
\usepackage{tikz}
\usepackage{tikz-cd}
\usetikzlibrary{matrix,calc}
\usepackage{hyperref}
\hypersetup{colorlinks,bookmarks=true,linktocpage=true,citecolor=blue,linkcolor=magenta}

\usepackage{eulervm}   

\allowdisplaybreaks 
\usepackage{color}
\usepackage[all]{xy} 

\usepackage{mathtools}
\usepackage{MnSymbol} 

\theoremstyle{plain}
\newtheorem{thm}{Theorem}[section]

\theoremstyle{definition}

\newtheorem{rem}[thm]{Remark}

\usepackage{enumerate}
\usepackage{enumitem}

\newcommand{\bRn}{\mathbb{R}^n}

\newcommand{\pd}{\partial}

\newcommand{\bP}{{\mathbb P}}

\newcommand{\bR}{{\mathbb R}}
\newcommand{\bS}{{\mathbb S}}

\newcommand{\cH}{{\mathcal H}}
\newcommand{\cI}{{\mathcal I}}

\newcommand{\cL}{{\mathcal L}}

\newcommand{\fR}{{\mathbf R}}

\newcommand{\wt}{\widetilde}

\def\barint_#1{\mathchoice
            {\mathop{\vrule width 6pt
height 3 pt depth -2.5pt
                    \kern -9.5pt
\intop \kern -4pt}\nolimits_{#1}}%
            {\mathop{\vrule width 5pt height
3 pt depth -2.6pt
                    \kern -6.5pt
\intop \kern -4pt}\nolimits_{#1}}%
            {\mathop{\vrule width 5pt height
3 pt depth -2.6pt
                    \kern -6pt
\intop \kern -4pt}\nolimits_{#1}}%
            {\mathop{\vrule width 5pt height
3 pt depth -2.6pt
          \kern -6pt \intop \kern -4pt}\nolimits_{#1}}}
          
           \def\bariint_#1{\mathchoice
            {\mathop{\vrule width 15pt
height 3 pt depth -2.5pt
                    \kern -15.8pt
\intop \kern -8pt\intop \kern -4pt}\nolimits_{#1}}%
            {\mathop{\vrule width 9pt height
3 pt depth -2.6pt
                    \kern -10.5pt
\intop \kern -8pt\intop \kern -4pt}\nolimits_{#1}}%
            {\mathop{\vrule width 9pt height
3 pt depth -2.6pt
                    \kern -10pt
\intop \kern -8pt\intop \kern -4pt}\nolimits_{#1}}%
            {\mathop{\vrule width 9pt height
3 pt depth -2.6pt
          \kern -8pt \intop \kern -10pt\intop \kern -4pt}
      \nolimits_{  #1}}}

\def\barintlim_#1{\mathchoice
            {\mathop{\vrule width 6pt
height 3 pt depth -2.5pt
                    \kern -8.8pt
\intop \kern -4pt}\limits_{#1}}%
            {\mathop{\vrule width 5pt height
3 pt depth -2.6pt
                    \kern -6.5pt
\intop \kern -4pt}\limits_{#1}}%
            {\mathop{\vrule width 5pt height
3 pt depth -2.6pt
                    \kern -6pt
\intop \kern -4pt}\limits_{#1}}%
            {\mathop{\vrule width 5pt height
3 pt depth -2.6pt
          \kern -6pt \intop \kern -4pt}\limits_{#1}}}
          
           \def\bariintlim_#1{\mathchoice
            {\mathop{\vrule width 15pt
height 3 pt depth -2.5pt
                    \kern -15.8pt
\intop \kern -8pt\intop \kern -4pt}\limits_{#1}}%
            {\mathop{\vrule width 9pt height
3 pt depth -2.6pt
                    \kern -10.5pt
\intop \kern -8pt\intop \kern -4pt}\limits_{#1}}%
            {\mathop{\vrule width 9pt height
3 pt depth -2.6pt
                    \kern -10pt
\intop \kern -8pt\intop \kern -4pt}\limits_{#1}}%
            {\mathop{\vrule width 9pt height
3 pt depth -2.6pt
          \kern -8pt \intop \kern -10pt\intop \kern -4pt}
      \limits_{  #1}}}
          
\renewcommand{\iint}{\int \kern -8pt\int}       

\newcommand{\RE}{\text{{\rm Re}}\,}

\numberwithin{equation}{section}
\makeatletter
\@namedef{subjclassname@2020}{\textup{2020} Mathematics Subject Classification}
\makeatother

\title{Evans-Lewis inequalities in dimension three}

\author{Yi C. Huang} 
\address{School of Mathematical Sciences, Nanjing Normal University, Nanjing 210023, People's Republic of China}
\email{Yi.Huang.Analysis@gmail.com}
\urladdr{https://orcid.org/0000-0002-1297-7674}

\author{Li Liu} 
\address{School of Mathematical Sciences, Yangzhou University, Yangzhou 225002, People's Republic of China}
\email{lliu@yzu.edu.cn}

\date{\today} 

\thanks{Research of the authors is supported by the National NSF grant of China (no. 11801274).
The first author (YCH) thanks Professor Tohru Ozawa (Waseda University) for kind support.}

\begin{document}

\begin{abstract}
In proving Rellich inequalities in the framework of equalities,
N. Bez, S. Machihara, and T. Ozawa obtained some interesting norm inequalities in the spirit of Evans and Lewis that compare the standard Laplacian with its radial and spherical components.
In this paper we give a simple unified proof and a strict improvement of these Evans-Lewis inequalities in the subtle dimension three case.
Our approach is robust and explains clearly the occurrence of the sharp constant.
\end{abstract}

\maketitle


\section{Introduction}

Denote by $|\cdot|$, $\|\cdot\|_2$ and $\langle\cdot,\cdot\rangle$ the usual notations for the Euclidean distance in $\bRn$, the $L^2$-norm on $\bRn$ and the inner product for $L^2(\bRn)$.
Recall the Rellich inequality:
\begin{equation} \label{e:Rellich}
\|\Delta f\|_2^2\geq\fR_n^2\left\|\frac{f}{|x|^2}\right\|_2^2,
\end{equation}
where $\fR_n=\frac{n(n-4)}{4}$.
See Rellich \cite{Rel56} and \cite{DavHin98, Yaf99, AdiGroSan06, TerZog07, KomOza09, GhoMor11, CalMus12, Avk18, Caz20, DuyLamPhi21, KelPinPog21, MetNegSobSpi21, BerGanRoy22} 
for various extensions in various contexts.
In particular, \eqref{e:Rellich} for $n=4$ becomes trivial.
If we wish to allow general functions in the Sobolev space $H^2(\bRn)$,
it is also necessary to restrict \eqref{e:Rellich} to $n\geq5$. 
For the Laplace operator $\Delta$, we have the standard decomposition
$$\Delta=\Delta_r+\frac{1}{|x|^2}\Delta_{\bS^{n-1}}=:\Delta_r+\Delta_s,$$
where $\Delta_r=\pd_r^2+\frac{n-1}{|x|}\pd_r$, 
with $\pd_r=\frac{x}{|x|}\cdot\nabla$ being the radial derivative,
and $\Delta_{\bS^{n-1}}$ denotes the Laplace-Beltrami operator on $\bS^{n-1}$.
Introduce the spherical derivatives
$$L_j=\pd_j-\frac{x_j}{|x|}\pd_r,\quad j=1, \cdots, n.$$
Thus we recognize the following relation
$$\frac{1}{|x|^2}\Delta_{\bS^{n-1}}=\sum_{j=1}^nL_j^2.$$
Note that $-\sum_{j=1}^nL_j^2$ is non-negative on $L^2(\bRn)$.

Recently in \cite{MacOzaWad17, BezMacOza23}, 
Machihara-Ozawa-Wadade and Bez-Machihara-Ozawa took the unusual perspective of understanding \eqref{e:Rellich} in the framework of equalities,
and found the following elegant results: 
for $n\geq2$ and $f\in C_0^\infty(\bRn\backslash\{0\})$,
\begin{equation} \label{e:equal}
\|\Delta f\|_2^2=\left\|\Delta_r f\right\|_2^2+
\left\|\sum_{j=1}^nL_j^2 f\right\|_2^2+2\fR_n\sum_{j=1}^n\left\|\frac{L_jf}{|x|}\right\|_2^2+2\left\langle-\sum_{j=1}^nL_j^2f_*,f_*\right\rangle,
\end{equation}
where $$f_*=\pd_rf+\frac{n-4}{2}\frac{f}{|x|}.$$
The last two terms in \eqref{e:equal} are non-negative for $n\geq4$, i.e.,
$$\RE\langle \Delta_r f, \Delta_sf\rangle\geq0,$$
or equivalently,
\begin{equation} \label{e:Hilbert}
\|\Delta f\|^2_2\geq\|\Delta_r f\|^2_2+\|\Delta_s f\|^2_2.
\end{equation}
These inequalities were first noticed by Evans and Lewis in \cite{EvaLew05}.
See \cite{Hua23} for a simpler proof not relying on the Bez-Machihara-Ozawa-Wadade equalities \eqref{e:equal}.

Unexpectedly, the Evans-Lewis inequalities \eqref{e:Hilbert} for $n=3$ are very subtle.

\begin{thm}[Bez-Machihara-Ozawa] \label{thm:BMO}
For $n=3$ and $f\in C_0^\infty(\mathbb{R}^3\backslash\{0\})$, we have
\begin{equation} \label{e:BMO1}
\|\Delta f\|^2_2\geq\|\Delta_r f\|^2_2
\end{equation}
and
\begin{equation} \label{e:BMO2}
\|\Delta f\|^2_2\geq\left(\frac{5}{8}\right)^2\|\Delta_s f\|^2_2.
\end{equation}
\end{thm}

Both \eqref{e:BMO1} and \eqref{e:BMO2} are sharp.
Our aim here is to give a simple unified proof of Theorem \ref{thm:BMO}, 
explaining the occurrence of the constant $\left(\frac{5}{8}\right)^2$ in a visible way.
The byproducts of our machinery include a strict improvement of \eqref{e:BMO1} 
and a weighted version of \eqref{e:Hilbert} for $n=3$ that interpolates \eqref{e:BMO2} and the improved form of \eqref{e:BMO1}. 

\begin{thm} \label{thm:BMO'}
For $n=3$ and $f\in C_0^\infty(\mathbb{R}^3\backslash\{0\})$, we have
\begin{equation} \label{e:BMO1'}
\|\Delta f\|^2_2\geq\|\Delta_r f\|^2_2+\frac14\|\Delta_s f\|^2_2.
\end{equation}
More generally, for $k_1\in[0,1]$, $k_2\in[0,25/64]$, and $9k_1+64k_2\leq25$, we have
\begin{equation} \label{e:BMO2'}
\|\Delta f\|^2_2\geq k_1\|\Delta_r f\|^2_2+k_2\|\Delta_s f\|^2_2.
\end{equation}
\end{thm}

\begin{rem}
It is desirable to extend Theorems \ref{thm:BMO}-\ref{thm:BMO'} to an abstract framework, say, on homogeneous groups 
as in Ruzhansky and Suragan \cite{RuzSur17, RuzSur19}.
The latter monograph also contains an excellent account of Hardy and Rellich inequalities.
We choose not to bother ourselves in this short paper to review such a huge literature.
\end{rem}

\section{Proof of Theorem \ref{thm:BMO'}}

Recall the following spectral computation of norms (\cite{BezMacOza23}): for $f\in C_0^\infty(\mathbb{R}^3\backslash\{0\})$,
\begin{equation} \label{e:sc1}
\left\|\sum_{j=1}^3L_j^2 f\right\|_2^2=\sum_{k=1}^\infty\mu_k^2\left\|\frac{\bP_kf}{|x|^2}\right\|_2^2,
\end{equation}

\begin{equation} \label{e:sc2}
\sum_{j=1}^3\left\|\frac{L_j f}{|x|}\right\|_2^2=\sum_{k=1}^\infty\mu_k\left\|\frac{\bP_kf}{|x|^2}\right\|_2^2,
\end{equation}

\begin{equation} \label{e:sc3}
\left\|\frac{f}{|x|^2}\right\|_2^2=\sum_{k=0}^\infty\left\|\frac{\bP_kf}{|x|^2}\right\|_2^2.
\end{equation}
Here $\bP_k$ is the orthogonal projection onto the closed subspace $\cH_k(\mathbb{R}^3)$ of $L^2(\mathbb{R}^3)$ spanned by spherical harmonics of order $k$
multiplied by radial functions, and
$$-\Delta_{\bS^{2}}\bP_k=\mu_k\bP_k,\quad \mu_k=k(k+1),\,k\geq0.$$
See the book by Dai and Xu \cite{DaiXu13} for details on spherical harmonics.

\textit{Proof of \eqref{e:BMO1'}}. According to \eqref{e:equal} for $n=3$, we have
$$\begin{aligned}
\|\Delta f\|_2^2&=\left\|\Delta_r f\right\|_2^2+(1-\alpha)
\left\|\sum_{j=1}^3L_j^2 f\right\|_2^2\\
&\qquad+\alpha\left\|\sum_{j=1}^3L_j^2 f\right\|_2^2-\frac32\sum_{j=1}^3\left\|\frac{L_jf}{|x|}\right\|_2^2\\
&\qquad+2\left\langle-\sum_{j=1}^3L_j^2f_*,f_*\right\rangle.
\end{aligned}$$
Using \eqref{e:sc1} and \eqref{e:sc2}, we have
$$\begin{aligned}
N_1(\alpha):&=\alpha\left\|\sum_{j=1}^3L_j^2 f\right\|_2^2-\frac32\sum_{j=1}^3\left\|\frac{L_jf}{|x|}\right\|_2^2\\
&=\sum_{k=1}^\infty\mu_k\left(\alpha\mu_k-\frac32\right)\left\|\frac{\bP_kf}{|x|^2}\right\|_2^2.
\end{aligned}$$
Since $\mu_k$ is increasing in $k$,
to ensure the non-negativity of $N_1(\alpha)$, it suffices to take $\alpha\geq\frac32/\mu_1=3/4$.
So, $1-\alpha\leq1/4$. Thus \eqref{e:BMO1'} follows by the non-negativity of $-\sum_{j=1}^3L_j^2$.

\textit{A simpler proof of \eqref{e:BMO2}}. 
We shall use the radial Rellich inequality for $n=3$
\begin{equation} \label{e:rad}
\left\|\Delta_r f\right\|_2^2\geq\frac{9}{16}\left\|\frac{f}{|x|^2}\right\|_2^2,
\end{equation}
see \cite[Theorem 3]{BezMacOza23} or \cite{HuaShi23}.
Similar to the proof of \eqref{e:BMO1'} and using \eqref{e:rad},
 $$\begin{aligned}
\|\Delta f\|_2^2&\geq(1-\alpha)
\left\|\sum_{j=1}^3L_j^2 f\right\|_2^2\\
&\qquad+\alpha\left\|\sum_{j=1}^3L_j^2 f\right\|_2^2-\frac32\sum_{j=1}^3\left\|\frac{L_jf}{|x|}\right\|_2^2+\frac{9}{16}\left\|\frac{f}{|x|^2}\right\|_2^2\\
&\qquad+2\left\langle-\sum_{j=1}^3L_j^2f_*,f_*\right\rangle.
\end{aligned}$$
Using \eqref{e:sc1}, \eqref{e:sc2} and \eqref{e:sc3} (dropping the term involving $\bP_0$), we have
$$\begin{aligned}
N_2(\alpha):&=\alpha\left\|\sum_{j=1}^3L_j^2 f\right\|_2^2-\frac32\sum_{j=1}^3\left\|\frac{L_jf}{|x|}\right\|_2^2+\frac{9}{16}\left\|\frac{f}{|x|^2}\right\|_2^2\\
&\geq\sum_{k=1}^\infty\left[\mu_k\left(\alpha\mu_k-\frac32\right)+\frac{9}{16}\right]\left\|\frac{\bP_kf}{|x|^2}\right\|_2^2.
\end{aligned}$$
Reasoning as before, to ensure $\mu_k\left(\alpha\mu_k-\frac32\right)+\frac{9}{16}\geq0$ for $k=1$, it suffices to take $\alpha\geq39/64$.
So, $1-\alpha\leq25/64$. Thus \eqref{e:BMO2} follows by the non-negativity of $-\sum_{j=1}^3L_j^2$.

\textit{Proof of \eqref{e:BMO2'}}. 
A tedious extension of above proof of \eqref{e:BMO2}. We omit the details.

\section{Evans-Lewis inequalities for generalized Laplacians}

Our machinery also works for the abstract setting (for our purpose, we focus\footnotemark
\footnotetext{It extends \textit{mutatis mutandis} to $\bRn$ with $n\geq2$.} on $\bR^3$) as considered in \cite{EvaLew05, BezMacOza23} for the following generalized Laplacian
$$\cL=\Delta_r+\frac{1}{r^2}\Lambda=:\Delta_r+\wt\Delta_s,$$
where $-\Lambda$ is a non-negative self-adjoint operator on $L^2(\bS^{2})$ whose spectrum is purely discrete with isolated eigenvalues
$\{\lambda_k\}_{k\in\cI}$ which may accumulate only at $+\infty$.
We shall use the same notation $\bP_k$ as before to denote the orthogonal projection onto the closed subspace of $L^2(\mathbb{R}^3)$ 
spanned by the $\lambda_k$-eigenspace of $-\Lambda$ multiplied by radial functions, 
thereby we have the following useful spectral relation
$$-\Lambda\bP_k=\lambda_k\bP_k.$$
We also assume that $0$ is an eigenvalue of $-\Lambda$ and write $\lambda_0=0$.

Using the spectral decomposition and dropping the term involving $\bP_0$, we have
$$\begin{aligned}
\|\cL f\|_2^2&=\left\|\Delta_r f\right\|_2^2+(1-\alpha)
\left\|\wt\Delta_s f\right\|_2^2\\
&\qquad+\alpha\left\|\wt\Delta_s f\right\|_2^2+\sum_{k\in\cI}2\RE\left\langle \Delta_rf,\wt\Delta_sf\right\rangle\\
&=(1-\alpha)
\left\|\wt\Delta_s f\right\|_2^2\\
&\qquad+\sum_{k\in\cI}\left(\left\|\Delta_r \bP_kf\right\|_2^2+\alpha\left\|\wt\Delta_s \bP_kf\right\|_2^2+2\RE\left\langle \Delta_r\bP_kf,\wt\Delta_s\bP_kf\right\rangle\right)\\
&\geq(1-\alpha)
\left\|\wt\Delta_s f\right\|_2^2\\
&\qquad+\sum_{k\in\cI\backslash\{0\}}\left(\fR_3^2\left\|\frac{\bP_kf}{|x|^2}\right\|_2^2+\alpha\lambda_k^2\left\|\frac{\bP_kf}{|x|^2}\right\|_2^2+2\fR_3\lambda_k\left\|\frac{\bP_kf}{|x|^2}\right\|_2^2\right).
\end{aligned}$$
In the last step, we used the Rellich inequality \eqref{e:rad} and \cite[Equation (2.15)]{EvaLew05}\footnotemark\footnotetext{This is a simple consequence of integration by parts.}.
Since $\lambda_k\neq0$ for $k\in\cI\backslash\{0\}$,
to ensure $\fR_3^2+\alpha\lambda_k^2+2\fR_3\lambda_k\geq0$, 
it suffices to take $$\alpha\geq\max_{k\in\cI\backslash\{0\}}-\fR_3^2\left(\frac{1}{\lambda_k}+\frac{1}{\fR_3}\right)^2+1=:-\fR_3^2\left(\frac{1}{\lambda_{k_0}}+\frac{1}{\fR_3}\right)^2+1,$$
which is equivalent to 
$$\frac{1}{1-\alpha}\geq \left(\frac{\lambda_{k_0}}{\lambda_{k_0}+\fR_3}\right)^2.$$ 
This leads to the following Evans-Lewis inequalities
$$\left\|\frac{\Lambda f}{|x|^2}\right\|_2^2\leq \left(\frac{\lambda_{k_0}}{\lambda_{k_0}+\fR_3}\right)^2\|\cL f\|^2_2.$$
which agree with \cite[Theorem 7]{BezMacOza23}.

\begin{rem}
The above simple approach is robust since we did not use the Bez-Machihara-Ozawa-Wadade equalities \eqref{e:equal}.
It shares some similarity\footnotemark\footnotetext{More precisely, both proofs use \cite[Equation (2.15)]{EvaLew05} instead of \eqref{e:equal}.} with the proof in \cite[Section 4.1]{BezMacOza23},
but the occurrence of the sharp constant is visible here.
\end{rem}

\bigskip

\section*{\textbf{Compliance with Ethical Standards}}

\bigskip

\textbf{Conflict of interest} The authors have no known competing financial interests 
or personal relationships that could have appeared to influence this reported work.

\bigskip

\textbf{Availability of data and material} Not applicable.

\bigskip

\bibliographystyle{alpha}
 
\bibliography{HuaY-LiuL-EvansLewisDimThree}

\end{document}